\documentclass[a4paper,12pt]{article}
\usepackage{float}
\usepackage{setspace}
\usepackage{graphicx}
\usepackage{amsfonts}
\usepackage{amssymb}
\usepackage{amsmath}
\usepackage{bbm}
\usepackage{wasysym}
\usepackage{multicol}
\usepackage{esint}
\usepackage[nottoc]{tocbibind}
\usepackage{array}
\usepackage{indentfirst}
\usepackage{authblk}
\usepackage{blindtext}
\usepackage[hmargin={2.5cm},vmargin={2.5cm}]{geometry}
\usepackage{scalefnt}
\usepackage{colortbl}
\usepackage[figurename=Fig.,
            font=scriptsize,
            hang,
            labelfont=bf,
            up,
            textfont=it,
            up]{caption} 
\usepackage{booktabs} 
\usepackage{lineno,hyperref}
\modulolinenumbers[5]
\definecolor{blue}{RGB}{00,204,204}
\definecolor{dgreen}{RGB}{51,102,00}
\hypersetup{colorlinks=true, citecolor=blue, linkcolor=blue,  urlcolor=blue}
\newtheorem{theo}{Theorem}[section]
\newtheorem{coro}[theo]{Corolary}

\newtheorem{prop}[theo]{Proposition}
\newcommand{\dem}{\no \textsc{Proof:~}}

\newtheorem{defi}{Definition}[section]

\newtheorem{rmk}{Remark}[section]

\newcommand{\cqd}{\nopagebreak\hfill\fbox{}}

\newcommand{\no}{\noindent}

\newcommand{\C}{\mathbb{C}}

\title{Numerical solutions of ordinary differential equations using Spline-Integral Operator}
\author[1]{Gustavo H. O. Salgado \thanks{gusalg@unifei.edu.br}~}
\author[1]{Jo\~ao P. R. Romanelli \thanks{joaoromanelli@unifei.edu.br}}
\affil[1]{Universidade Federal de Itajub\'{a}, Campus Itabira.
Rua Irm\~{a} Ivone Drumond, 200.}

\date{\today}

\makeindex

\begin{document}
\maketitle
\abstract{
In this work, we introduce a novel numerical method for solving initial value problems associated with a given differential. Our approach utilizes a spline approximation of the theoretical solution alongside the integral formulation of the analytical solution. Furthermore, we offer a rigorous proof of the method's order and provide a comprehensive stability analysis. Additionally, we showcase the effectiveness method through some examples, comparing with Taylor's methods of same order.
}   

\section{Introduction}

An Initial Value Problem (IVP) given by 
\begin{equation}\label{eq:ivp}
    \left\{
    \begin{array}{ccl}
    y'(t) & = & f(t,y(t)), \\
    y(t_{0}) & = & y_{0},
    \end{array}
    \right.
\end{equation}
consists of find a function $y(t)$ that solves the Ordinary Differential Equation (ODE) in IVP \eqref{eq:ivp} such that $y(t_0) = y_{0}$.
It is well known that if $f(\cdot, y)$ is Lipschitz continuous in some rectangle around the initial condition $(t_0,y_0)$, there is a unique solution $y$ of Eq.\,\eqref{eq:ivp} \cite{coddington_1955}. 
In most situations, finding the explicit solution of a given IVP can be awkward or even impossible.
In order to address this difficult, numerical methods are used and play an important role in most applications involving differential equations.
Since numerical methods for solving an ODE generally result in approximations of the theoretical solution, the study and analysis of errors play another important role in these methods.
An important class of numerical methods widely used to approximate solutions of IVP \eqref{eq:ivp} is the Runge--Kutta (RK) methods \cite{runge_1895,kutta_1901}.
They are obtained using Taylor's expansion of a theoretical solution $y(t)$ of \eqref{eq:ivp} written in terms of the vector field $f(t,y)$ and comparing with an incremental function that yields a nonlinear algebraic system with infinitely many solutions, where each solution provides a set of parameters for the incremental function.
A general methodology based on the use of tables and tree diagrams was developed by Butcher (1963) \cite{butcher_1963}, which simplifies the algebraic operations involved in Taylor's expansions to obtain the RK parameters. 
%
%

%
Another important class of numerical methods to solve ODEs is the linear multi-step methods \cite{quarteroni_2006}. 
These methods are divided into explicit and implicit categories. 
In an implicit multi-step, the numerical solution is obtained by solving a nonlinear equation that depends on $f(t,y)$.
This solution can be found by employing a Newton method to solve the nonlinear equation.
To avoid the evaluation of $f(t,y)$ at each step where the Newton method is applied, a prediction is obtained from an explicit multi-step method followed by a correction using a fixed-point scheme, leading to Predictor-Corrector Methods.
However, the application of multi-step methods relies on knowing a number of solutions of the IVP, which are generally obtained from the initial condition $y(t_0)=y_0$ using other methods, such as  Taylor or RK methods, which are one-step methods \cite{butcher_2016}. 
Butcher (2000) \cite{butcher_2000} provides an overview of some major developments of Runge-Kutta and multi-step procedures over the last century. 
The spline approach has been used to obtain numerical solutions for IVPs. Typically, the basic idea is consider a spline as an approximation of the solution.
The spline depends on a parameter associated with the term of the larger degree, which is computed. Loscalzo (1969) \cite{loscalzo_1963} used such an approach to obtain numerical approximations of first-order ODEs. In particular, it is demonstrated how to obtain the trapezoidal rule from spline of degree two and the Milne-Simpson method from a spline of degree three. 
In \cite{aziz_1981}, approximations for second order value problems and boundary problems for periodic solutions are obtained.
The spline approach is also used to approximate solution of nonlinear matrix ordinary differential problems \cite{defez_2021,tung_2022}.
A spline collocation method is proposed in \cite{mazzia_2006a}, where the collocation points coincide with the knots in  the mesh. The method, called B-Spline (BS) Method, compute the spline values at the knots and is associated a class of linear multi-step methods. 
The BS Method can also be applied for a non-uniform meshes \cite{mazzia_2006b}, and classes of quadrature formulas associated to the BS Method is studied in \cite{mazzia_2012}. 
In \cite{falini_2023}, the Spline Quasi-Interpolation method, based on B-Splines, is presented and is applied to numerical solution of Gauss-Lobatto and Gauss-Legendre methods using Runge-Kutta methods.
In this work, we propose a one-step method derived from a function obtained through spline approximation of the solution of the IVP combined with its integral formulation. We define a spline that approximates the solution, but depends on an unknown parameter, which will be computed by a fixed-point iteration and taken as the approximated solution of the IVP. 
We also demonstrate that there is an interval containing the initial condition where the approximate solution given by the proposed method is the unique solution in that interval.
Moreover, the method has order $(m+1)$ when the spline approximation of the solution is expressed using derivatives up to order $(m-1)$. Additionally, we present numerical experiments where the approximation obtained does not depend on any numerical integral method. 
This work is organized as follow: in Section \ref{sec:Spline_Int_form}, we define the Spline Intregal Operator and prove its convergence; in Section \ref{sec:SIO_order}, we discuss the order of the SIO method; in Section \ref{sec:stab_analysis}, we analyze the stability; in Section \ref{sec:num_examples}, we present some numerical experiments; and finally, in Section \ref{sec:conclusions}, we provide the conclusion and final remarks.
%
%
\section{The Spline-Integral Operator}\label{sec:Spline_Int_form}
Let $y(t)$ the theoretical solution of IVP \eqref{eq:ivp}.
If $y(t)$ has continuous derivatives up to $m\!-\!1$ at $t=t_0$, there is a unique spline of degree $m$ given by 
\begin{equation}\label{eq:spline}
S_m(t_0,t,w) = \sum_{k=0}^{m-1} \dfrac{1}{k!}\dfrac{d^ky(t_0)}{dt^k}\,(t-t_0)^k + s_m(t-t_0)^m,
\end{equation}
where $s_m$ satisfies $S_m(t_0,t_0+h,\,w)=w$ for a given pair $h$ and $w$.
Is straightforward that 
\begin{equation}\label{eq:sm}
 s_m = \dfrac{1}{h^m}\left(w-\sum_{k=0}^{m-1} \dfrac{1}{k!}\dfrac{d^ky(t_0)}{dt^k}\,h^k\right).
\end{equation}
%
\begin{rmk}
By the application of the chain rule, if $f(t,y)$ has continuous partial derivatives up to order $m\!-\!1$ at $(t_0,y_0)$, then the spline $S_m(t_0,t,w)$ can be obtained from IVP \eqref{eq:ivp}.  
\end{rmk}\medskip
%
\begin{rmk}
The Taylor polynomial of degree $m-1$ of a function $y(t)$ at $t=t_0$ is given by
\begin{equation}\label{eq:taylor_m-1}
 T_{m-1}(t_0,t) = \sum_{k=0}^{m-1} \dfrac{1}{k!}\dfrac{d^ky(t_0)}{dt^k}\,(t-t_0)^k = y(t_0) + \sum_{k=1}^{m-1} \dfrac{d^{k-1}f(t_0,y(t_0))}{dt^{k-1}}(t-t_0)^{k}.
\end{equation}
Hence, $s_m=(w-T_{m-1}(t_0,t_0+h))/h^m$ and $$S_m(t_0,t,w)=T_{m-1}(t_0,t)+s_m(t-t_0)^m.$$
\end{rmk}\medskip
%

We will define an integral operator inspired on the integral form of the solution of IVP \eqref{eq:ivp}, 
\begin{equation}\label{eq:int_form_I}
 y(t) = y(t_0) + \int_{t_0}^{t}\,f(\tau,y(\tau))\,d\tau,
\end{equation}
where $y(\tau)$ will be replaced by the spline defined in Eqs.\,\eqref{eq:spline} and \eqref{eq:sm}.
However, we need to establish conditions for $h$ and $w$ such that $f(t,\,S_m(t_0,t,w))$ is well-defined. 
In the following proposition, we will demonstrate conditions under which $f(t,S_m(t_0,t,w))$ is defined.
%
\begin{prop}\label{prop:SsubsetR}
Let the rectangle 
\begin{equation}\label{eq:R_TEU}
 \mathcal{R} = \{(t, y); |t-t_0| \leq a, |y-y_0| \leq b\},
\end{equation}
where $f(t,y)$ is continuous. 
If the theoretical solution $y(t)$ of IVP \eqref{eq:ivp} has $m-1$ continuous derivatives at $t=t_0$, then there exist $h$ and $w$ such that $(t,S_m(t_0,t,w))\subset\mathcal{R}$ for all $|t-t_0|\leq |h|$.
\end{prop}
%
\dem Since Taylor polynomial $T_{m-1}(t_0,t)$ is continuous at $t=t_0$, there is a $\delta>0$ such that
\begin{equation*}
 |T_{m-1}(t_0,t)-y_0|<\dfrac{b}{3} 
\end{equation*}
for each $|t-t_0|<\delta$.
For a given $h$ and $w$, such that $0<|h|<\delta$ and $|w-y_0|\leq \dfrac{b}{3},$ we have
\begin{eqnarray*}
 |S_m(t_0,t,w)-y_0|   & = &  |T_{m-1}(t_0,t) + s_m(t-t_0)^m-y_0| \\[0,3cm]
                    & \leq & |T_{m-1}(t_0,t) -y_0| + |s_m(t-t_0)^m| \\[0,3cm]
                    & \leq & |T_{m-1}(t_0,t) -y_0| + |s_m h^m| \\[0,3cm]
                      & = & |T_{m-1}(t_0,t) -y_0| + |w - T_{m-1}(t_0,t_0+h)| \\[0,3cm]
                      & = & |T_{m-1}(t_0,t) -y_0| + |w - y_0| + |T_{m-1}(t_0,t_0+h)-y_0| \\[0,3cm]
                      & = & \dfrac{b}{3}+\dfrac{b}{3}+\dfrac{b}{3}=b,
\end{eqnarray*}
for all $|t-t_0|\leq|h|$.
\cqd\medskip
%

%
\begin{rmk}
In Proposition \ref{prop:SsubsetR}, the existence of $h$ is guaranteed by the continuity of spline $S_m(t_0,t,w)$ for any $w$. However, in the proof of the aforementioned proposition, we show that $h$ can be chosen to account for the Taylor polynomial used in the spline definition. Also, for an appropriated choice of the independent values $w$ and $h$, the graph of $S_m(t_0,t,w)$ for $|t-t_0|\leq|h|$ is entirely contained in $\mathcal{R}$.
\end{rmk}\bigskip
%
%

The spline $S_m(t_0,t,w)$ will be used to define the integral operator as follows.
%
\begin{defi}
For a given  IVP \eqref{eq:ivp}, where $f(t,y)$ is continuous in a rectangle $\mathcal{R}$, and for a fixed $h$, the Spline-Integral Operator (SIO) $G_h(w)$ is defined as 
 \begin{equation}\label{eq:G_operator}
 G_h(w) = \,y_0 + \int_{t_0}^{t_0+h}\,f(t,S_m(t_0,t,w))\,dt,
\end{equation}
where $h$ and $w$ are such that $(t,S_m(t_0,t,w))\subset\mathcal{R}$ for all $|t-t_0|\leq |h|$.
\end{defi}
%

In the next theorem, we shall see that if $f(t,y)$ is continuous and Lipschitz with respect the variable $y$ in $\mathcal{R}$, then there is a suficient small $h$ such that $G_h(w)$ has a unique fixed point.
%
\begin{theo}\label{th:G_contraction}
Let $f(t,y)$ continuous in $\mathcal{R}$ (as defined in Proposition \ref{prop:SsubsetR}) and Lipschitz with respect variable $y$ with constant $K$, and $h$ such that 
$$\left|T_{m-1}(t_0,t)-y_0\right|<\dfrac{b}{3}, ~ \mbox{for all}~ |t-t_0|\leq|h|.$$
If $|f|\leq M$ in $\mathcal{R}$, then 
%
$G_\alpha:I\longrightarrow I, ~~I=[y_0-b/3,y_0+b/3],$
%
%
has a unique fixed point whenever $|\alpha|<\min\left\{|h|,\,\dfrac{b}{3M},\,\dfrac{m+1}{K}\right\}$. 
\end{theo}
%
\dem For each fixed $|\alpha|<\min\left\{|h|,\dfrac{b}{3M},\dfrac{m+1}{K}\right\}$ we have
\begin{eqnarray*}
 |G_\alpha(w)-y_0| & = & \left|\int_{t_0}^{t_0+\alpha}\,f(t,S_m(t_0,t,w))\,dt\right| \leq M|\alpha| <  \dfrac{b}{3},
\end{eqnarray*}
for each $w\in I$.
In the sequel, we shall show that $G_\alpha(w)$ is a contraction.
Given $w_1,w_2\in I$, then
\begin{eqnarray*}
 |G_\alpha(w_1)-G_\alpha(w_2)| & = & \left|\int_{t_0}^{t_0+\alpha}[f(t,S_m(t_0,t,w_1))-f(t,S_m(t_0,t,w_2))]dt\right| \\[0,3cm]
                & \leq & \int_{t_0}^{t_0+\alpha}|f(t,S_m(t_0,t,w_1))-f(t_0,t,S_m(t_0,t,w_2))|dt \\[0,3cm] 
                & \leq & \int_{t_0}^{t_0+\alpha} K|S_m(t_0,t,w_1)-S_m(t_0,t,w_2)|dt \\[0,3cm] 
                & \leq & K\int_{t_0}^{t_0+\alpha} \left|\dfrac{(w_1-w_2)(t-t_0)^m}{\alpha^m}\right|dt \\[0,3cm] 
 %
%
                & \leq & \dfrac{K|\alpha|}{m+1}|w_1-w_2|. 
\end{eqnarray*}
Therefore, $G_\alpha(w)$ is a contraction since $|\alpha|<\dfrac{m+1}{K}$ and, from the Banach fixed point theorem \cite{banach_1922}, $G_\alpha(w)$ has a unique fixed point.
\cqd\medskip
%

\begin{rmk} It is well known that if $\partial f/ \partial y$ is continuous in $\mathcal{R}$, then $f(t,y)$ is Lipschitz. In that case, $K$ can be replaced by $\displaystyle\max\{\partial f/\partial y\}$ on $\mathcal{R}$ in Theorem \ref{th:G_contraction}.
\end{rmk}\medskip 
%

%
In the next section, we will show that the fixed point of $G_h(w)$, obtained as the limit of sequence $w_{n+1}=G_h(w_n)$, with $w_0=y_0$, is an approximation of the theoretical solution $y(t_0+h)$ of order $m+1$ when $h$ satisfies the hypothesis of Theorem \ref{th:G_contraction}.
Additionally, we will extend the method for a set of nodes in an interval containing $t_0$ as a method to solve IVPs.
%
%
%
%
%
\section{Order approximation of SIO}\label{sec:SIO_order}
The SIO given by Eq.\,\eqref{eq:G_operator} can be used as an approximation of the theoretical solution $y(t_0+h)$ of an IVP \eqref{eq:ivp}.
In the sequel, we shall show that each term of the sequence $w_{n+1}=G_h(w_n)$, with $n>0$ and $|w_0-y_0|\leq b/3$ is at least an  approximation of order $m$. 
%
\begin{prop}\label{prop:G_h(w)_approx_m}
 If $f(t,y)$, $w$, and $h$ satisfy the hypothesis of Theorem \ref{th:G_contraction} and the theoretical solution $y(t)$ has $m+1$ continuous derivatives, then
 $|G_h(w)-y(t_0+h)|$ is $\mathcal{O}(h^{m+1})$. 
\end{prop}
\dem Let $g_{{w}}(t)$ defined by
\begin{equation*}
 g_{{w}}(t) = y_0 + \int_{t_0}^{t}f(\tau,S_m(t_0,\tau,{w}))d\tau.
\end{equation*}
By definition of $S_m(t_0,t,{w})$, $\dfrac{d^kS_m(t_0,t,{w})}{dt^k}\Big|_{t=t_0}=\dfrac{d^ky(t_0)}{dt^k}$ for $k<m$, and $$\dfrac{d^mS_m(t_0,t,w)}{dt^m}\Big|_{t=t_0}=m!s_m,$$ by Eq.\,\eqref{eq:spline}. 
Then, 
\begin{equation}\label{eq:prop_G_h(w)_approx_m}
 \dfrac{d^kg_{{w}}(t)}{dt^k}\Big|_{t=t_0} = \dfrac{d^ky(t_0)}{dt^k}, \mbox{ for } k\leq m. 
\end{equation}
Also, 
\begin{equation*}
\dfrac{d^{m+1}g_{{w}}(t)}{dt^{m+1}}\Big|_{t=t_0} = \dfrac{d^{m+1}y(t_0)}{dt^{m+1}} + \dfrac{\partial f(t_0,y_0)}{\partial y}\left(m!s_m -\dfrac{d^my(t_0)}{dt^m}\right),
\end{equation*}
by application of $(m+1)$th derivative on the integral form of $y(t)$, given by Eq.\,\eqref{eq:int_form_I}, and comparing with $\dfrac{d^{m+1}g_w(t)}{dt^{t^{m+1}}}$ at $t=t_0$.
Then, taking the Taylor series of $g_{{w}}(t)$ and the Taylor series of theoretical solution $y(t)$ at $t=t_0$ up to order $m+1$, 
we have
\begin{eqnarray*}
 g_{{w}}(t_0+h)-y(t_0+h) & = & \dfrac{\partial f(t_0,y(t_0))}{\partial y}\left(m!s_m -\dfrac{d^my(t_0)}{dt^m}\right)\dfrac{h^{m+1}}{(m+1)!} + \mathcal{O}(h^{m+2}),
\end{eqnarray*}
which proves the assertion, since $G_h(w)=g_w(t_0+h)$ . 
\cqd\medskip 
%

%
The order of approximation at the fixed point of $G_h(w)$ can be improved as a consequence of the proof of Proposition \ref{prop:G_h(w)_approx_m}. \bigskip
%

%
\begin{coro}\label{coro:G_h(w)_approx_m}
If $\tilde{w}$ is the fixed point of $G_h(w)$, then $|\tilde{w}-y(t_0+h)|$ is $\mathcal{O}(h^{m+2})$.
\end{coro}
\dem Since $\tilde{w}=G_h(\tilde{w})$, then $\tilde{w}=g_{\tilde{w}}(t_0+h)$. From Eqs.\,\eqref{eq:sm} and \eqref{eq:prop_G_h(w)_approx_m} it follows that 
\begin{eqnarray*}
 m!s_m - \dfrac{d^my(t_0)}{dt^m} & = & \dfrac{m!}{h^m}\left(\tilde{w}-\sum_{k=0}^{m-1} \dfrac{1}{k!}\dfrac{d^ky(t_0)}{dt^k}\,h^k- \dfrac{d^my(t_0)}{dt^m}\dfrac{h^m}{m!}\right) \\[0,3cm]
                                 & = & \dfrac{m!}{h^m}\left(g_{\tilde{w}}(t_0+h)-\sum_{k=0}^{m} \dfrac{1}{k!}\dfrac{d^ky(t_0)}{dt^k}\,h^k\right) \\[0,3cm]
                                 &  = & \dfrac{m!}{h^m}\dfrac{d^{m+1}g_{\tilde{w}}(\xi)}{dt^{m+1}}\dfrac{h^{m+1}}{(m+1)!},
\end{eqnarray*}
for some $\xi\in[t_0,t_0+h]$. Hence,
\begin{eqnarray*}
 \tilde{w}-y(t_0+h) & = & \dfrac{\partial f(t_0,y(t_0))}{\partial y} \left(\dfrac{d^{m+1}g_{\tilde{w}}(\xi)}{dt^{m+1}}\dfrac{h}{m+1}\right)\dfrac{h^{m+1}}{(m+1)!} + \mathcal{O}(h^{m+2})\\[0,3cm]
                    & = & \dfrac{\partial f(t_0,y(t_0))}{\partial y}\dfrac{d^{m+1}g_{\tilde{w}}(\xi)}{dt^{m+1}}\dfrac{h^{m+2}}{(m+1)^2\,(m!)} + \mathcal{O}(h^{m+2}) ~=~  \mathcal{O}(h^{m+2}).
\end{eqnarray*}
\cqd\bigskip
%

%
Therefore, the approximation obtained from $G_h(w)$ of $y(t_0+h)$ is of order $m$, and it becomes an approximation of order $m+1$ for $\tilde{w}=G_h(\tilde{w})$. Henceforth, if the spline used in the SIO is of degree $m$, we will refer to it as \textbf{SIO of order} $\boldsymbol{m+1}$, denoted as \textbf{SIO}$\boldsymbol{(m+1)}$. 
The SIO of order $m+1$ can be extended to determine numerical solutions for an IVP \eqref{eq:ivp} over a discretization nodes $t_{i+1} = t_i + h$, with $i=0,\dots,N-1$,
using the implicit recursion schema
\begin{equation}\label{eq:G_opt_ti}
\left\{\begin{array}{l}
  \tilde{w}_0=y(t_0),\\[0,3cm]
  %
  \tilde{w}_{i+1}= \tilde{w}_{i} + \displaystyle\int_{t_i}^{t_i+h}f(t,S_m(t_i,t,\tilde{w}_{i+1}))\,dt,~i=0,\dots,N-1,
\end{array}\right.
\end{equation}
where the Taylor polynomial part of the spline $S_m(t_i,t,w)$ is evaluated by replacing $y(t_i)$ with $\tilde{w}_i$ in the right-hand side (RHS) of Eq.\,\eqref{eq:taylor_m-1} at each step $i$.
Clearly, $\tilde{w}_{i+1}$ is the unique fixed point of the SIO($m$+1) in the second equation of \eqref{eq:G_opt_ti}, if $t_i$, $\tilde{w}_i$, and $h$ satisfy an analogous condition to that of Theorem \ref{th:G_contraction}.
In the next section, we will provide a comprehensive stability analysis for this method, encompassing cases for each $m \geq 1$.
Additionally, we will present a detailed and explicit derivation of the SIO method outlined in Eq.\,\eqref{eq:G_opt_ti} for a linear homogeneous differential equation.
This derivation will serve as a foundation for demonstrating the stability analysis.
%
%
\section{Stability Analysis}\label{sec:stab_analysis}
The stability of the proposed method is determined by examining numerical solutions of the IVP
\begin{equation}\label{eq:linear_ivp}
 y'=\lambda y, ~~ y(t_0)=y_0,
\end{equation}
where $\lambda$ is a complex parameter. For a given $m$, the associated spline of Eq.\,\eqref{eq:linear_ivp} at $t_0$ is given by 
\begin{equation*}
 S_m(t_0,t,w) = y_0\,\sum_{k=0}^{m-1} \dfrac{\lambda^k(t-t_0)^k}{k!} + s_m(t-t_0)^m,
\end{equation*}
where $s_m=\dfrac{1}{h^m}\left(w-y_0\displaystyle\sum_{k=0}^{m-1}\dfrac{(\lambda h)^k}{k!}\right)$.
Then, using the schema in Eq.\,\eqref{eq:G_opt_ti} for $i=1$ ($\tilde{w}_0=y_0$), we obtain the Spline-Integral Operator
\begin{eqnarray*}
 G_h(w)  & = & y_0 + \int_{t_0}^{t_0+h} \lambda\left(y_0\,\sum_{k=0}^{m-1} \dfrac{\lambda^k(t-t_0)^k}{k!} + s_m(t-t_0)^m\right)dt \\[0,3cm]
         & = & y_0 + y_0\lambda\,\sum_{k=0}^{m-1} \dfrac{\lambda^k h^{k+1}}{(k+1)!} + \dfrac{s_m\lambda h^{m+1}}{m+1} \\[0,3cm]
         & = & y_0\left(1+\sum_{k=0}^{m-1} \dfrac{(\lambda h)^{k+1}}{(k+1)!}\right) + \dfrac{\lambda h}{m+1}\left(w-y_0\displaystyle\sum_{k=0}^{m-1}\dfrac{(\lambda h)^k}{k!}\right) \\[0,3cm]
         & = & \dfrac{\lambda h w}{m+1} + y_0\left(1+\sum_{k=0}^{m-1} \dfrac{(m-k)(\lambda h)^{k+1}}{(m+1)(k+1)!}\right). \\[0,3cm]
\end{eqnarray*}
The fixed point of $G_h(w)$ at $t_0$, \textit{i.e.}, for $i=1$, is given by
\begin{equation}
 \tilde{w}_1 = \dfrac{y_0}{m+1-\lambda h}\left(m+1+\sum_{k=0}^{m-1} \dfrac{(m-k)(\lambda h)^{k+1}}{(k+1)!}\right).
\end{equation}
By induction in $i$, we have
\begin{equation}
 \tilde{w}_i = \dfrac{y_0}{(m+1-\lambda h)^i}\left(m+1+\sum_{k=0}^{m-1} \dfrac{(m-k)(\lambda h)^{k+1}}{(k+1)!}\right)^i.
\end{equation}
This shows the following proposition
\begin{prop}\label{prop:stab_SIO}
The SIO($m$+1) is stable if, and only if,  
\begin{equation}\label{eq:stab_SIO}
\dfrac{1}{|m+1-\lambda h|}\left\vert m+1+\sum_{k=0}^{m-1} \dfrac{(m-k)(\lambda h)^{k+1}}{(k+1)!}\right\vert<1,~~ \lambda h\in\C.
\end{equation}
\cqd
\end{prop}\bigskip

In the following discussion, we will utilize Proposition \ref{prop:stab_SIO} to define the stability region of SIO($m+1$) and compare its stability with that of Taylor’s methods across various values of  $m$.

The stability of the SIO(2) is equivalent to that of the implicit trapezoidal method \cite{quarteroni_2006}. In fact, we have 
\begin{equation}
 \tilde{w}_i = y_0\,\left(\dfrac{2+\lambda h}{2-\lambda h}\right)^i, ~ \mbox{for} ~ m=1 ~ \mbox{and} ~ i\geq 0.
\end{equation}
Hence, the absolute stability region for SIO(2) is the half semi-plane $\{\lambda h\in\C; ~ \mathrm{Re}(\lambda h)<0\}$.

\begin{figure}[H]
 \centering
 \includegraphics[trim={0 0 3cm 0}, scale=1.2]{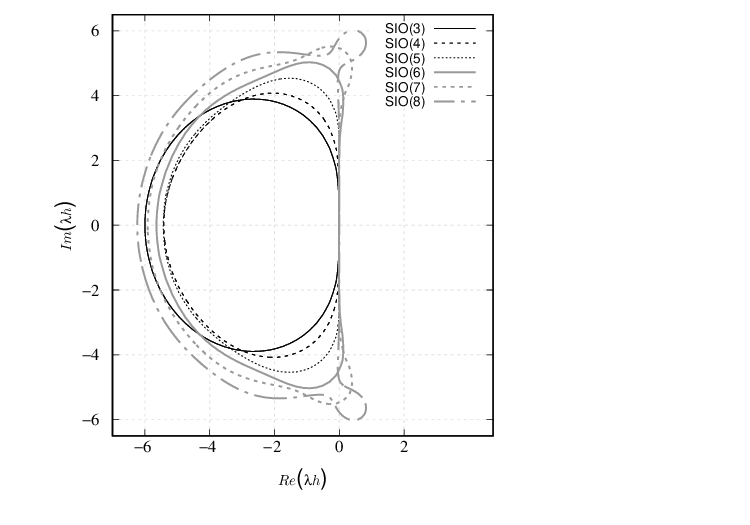}\vspace{-0,25cm}

 \caption{Stability region of SIO($m$+1) method for $m\in\{2,\dots,7\}$.}\label{fig:stab_SIO}
\end{figure}

In Figure \ref{fig:stab_SIO}, we illustrate the stability regions of the SIO($m$+1) for $m\in\{2,\dots,7\}$. As depicted, $m=7$ marks the first instance where the region encompasses the stability region for $m=2$.      

\begin{figure}[H]
 \centering
 \includegraphics[trim={0 0 3cm 0}, scale=1.2]{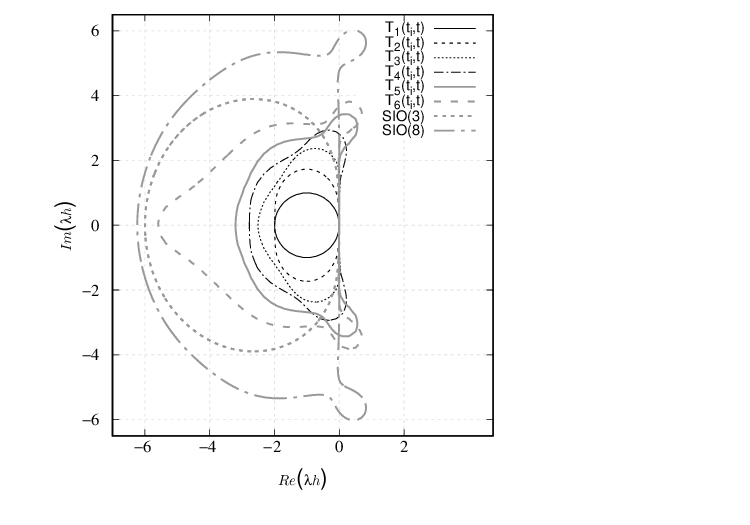}\vspace{-0,25cm}

 \parbox{10cm}{\caption{Stability regions of SIO(3) and SIO(8) compared with stability regions of Taylor's methods of order less than 6.}\label{fig:stab_SIO_Tay}}
\end{figure}

In Figure \ref{fig:stab_SIO_Tay}, we compare the stability regions of the SIO($m$+1) for $m=2$ and $m=7$ with Taylor methods of order $1$ through $6$, which correspond to $m=2,\dots,7$. 
It is observed that The stability region of SIO(3) exceeds those of Taylor's methods of order less than 6. Moreover, the stability region  of SIO(3) encompasses a significant portion of the stability regions defined by these Taylor methods.      
The stability analysis of the SIO($m+1$) method developed in this section was facilitaded by the linear nature of the problem given by Eq.\,\eqref{eq:linear_ivp}. This is due to the integral of $f(t, S(t_i, t, w))$ was obtained in terms of elementary functions, and the fixed point of $G_h(w)$ was explicitly written, leadind us to the explicit derivation of SIO($m+1$) for the IVP \eqref{eq:linear_ivp}.
In the next section, the proposed method described  by Eqs.\,\eqref{eq:G_opt_ti} will be reformulated as an iterative method. Here, the integral will be computed using a known numerical method, and the fixed point will be recursively determined. This approach uses the contraction property of the SIO($m$+1).
The proposed method will be demonstrated through some examples, and the obtained results will be contrasted with the Taylor's method of order $m+1$, $T_{m+1}(t_i,t)$, which shares the same order of the proposed method. 
Additionally, we will compare it with Taylor's method of order $m\!-\!1$, $T_{m-1}(t_i,t)$, which is employed whithin the $S_m(t_i,t,w)$ spline.  
%
%
%
\section{Numerical Experiments}\label{sec:num_examples}
In this section, we will present the results of numerical experiments obtained used  SIO method developed in Sections \ref{sec:Spline_Int_form} and \ref{sec:SIO_order}, solving problems for which the analytical is known.
The SIO($m$+1) can be recursively applied to generate a series of numerical approximations, $\tilde{w}_i$, which approximate the theoretical solution $y(t_i)$ of an  IVP  \eqref{eq:ivp} across discretization nodes $t_i=t_0+ih$, with $i=0,\dots,N$. This is achieved through the following procedure 
\begin{equation}\label{eq:num_SIO}
 \left\{\begin{array}{l}
   \tilde{w}_0 \,=\, y_0, \\[0,75cm]
   %
   %
   ~~~ {w}_{i+1}^{[0]} \,=\, \tilde{w}_i, \\[0,3cm]
   ~~~ {w}_{i+1}^{[\ell]} \,=\, \tilde{w}_i +\displaystyle\int_{t_i}^{t_i+h}\!f\left(t,S_m\left(t_i,t,{w}_{i+1}^{[\ell-1]}\right)\right)dt,\, \ell=1,\dots,M,\\[0,5cm]
  \tilde{w}_{i+1}=w_{i+1}^{[M]}, ~~ i=0,\dots,N-1,
\end{array}\right.
\end{equation}\medskip

\noindent for a sufficient large $M$ or when $\left|w_{i+1}^{[\ell]}-w_{i+1}^{[\ell+1]}\right|$ is less than a tolerance in each step $i$. 
%

%
\begin{rmk}
From Proposition \ref{prop:SsubsetR} and Theorem \ref{th:G_contraction}, if the function $f(t,y)$ has continuous partial derivatives up to order $m\!-\!1$ in a step $(t_i,\tilde{w}_i)$, then $f(t,y)$ is Lipschitz continuous in a rectangle
\begin{equation*}
\mathcal{R}_i=\{(t,y);~|t-t_i|\leq a,~ |y-\tilde{w}_i|\leq b\},
\end{equation*}
and there is an $h>0$ where the approximation $\tilde{w}_{i+1}$ can be obtained.
\end{rmk}
%

%
In each example below, the degree of the spline $S_m(t_i,t,w)$ will be set as $m=3$, meaning the SIO has an approximation of order $4$.
To demonstrate the efficiency and robustness of the proposed method, we will also present results obtained using Taylor's methods of orders 2 and 4.  
Furthermore, we will discuss the cases where the integral in Eq.\eqref{eq:num_SIO} can be explicitly obtained and instances where it is approximated using numerical methods.

%
\begin{rmk}
In a general scenario, the integral of $f(t,S(t_i,t,w))$ may not yield expressions in terms of elementary functions. Consequently, an explicit expression for the fixed point of $G_h(w)$ cannot always be obtained.
\end{rmk}\medskip
%

\subsection*{Example 1}
Let the IVP 
\begin{equation}\label{eq:ex01}
 \left\{\begin{array}{l}
 y'=-y+t+2, \\[0,3cm] y(0)=2.
\end{array}\right. 
\end{equation}\bigskip
 The unique solution of the IVP \eqref{eq:ex01} is $y(t)=(t+1)+e^{-t}$.
In this case, the integral in the procedure outlined by Eq.\,\eqref{eq:num_SIO} can be calculated analytically.
The errors between the analytical solution and the numerical solutions obtained using $T_2(t-h,t)$, $T_4(t-h,t)$, and  SIO(4) are presented in Table \ref{tab:example_01}.
%
\begin{table}[H]
\footnotesize
\centering
\parbox{12cm}{\caption{Absolute error between analytical and numerical solutions for IVP $y'=-y+t+2$, $y(0)=2$, using $T_2(t-h,t)$, $T_4(t-h,t)$, and  SIO(4) ($m=3$), with $h=0.1$.}\label{tab:example_01}}
\begin{tabular}{ccccc}
\toprule
\multicolumn{2}{c}{Solution} &\multicolumn{3}{c}{Error} \\
\cmidrule(r){1-2} \cmidrule(r){3-5}
  $t$ & $y(t)$  & $T_2(t-h,t)$ & $T_4(t-h,t)$ & SIO(4) \\
 \midrule
 $0.1$ &  $2.00483741803596$ & $1.63\times10^{-4}$  & $8.20\times10^{-8}$ & $1.97\times10^{-8}$   \\[0,15cm]
 $0.2$ &  $2.01873075307798$ & $2.94\times10^{-4}$  & $1.48\times10^{-7}$ & $3.56\times10^{-8}$   \\[0,15cm]
 $0.3$ &  $2.04081822068172$ & $3.99\times10^{-4}$  & $2.01\times10^{-7}$ & $4.83\times10^{-8}$   \\[0,15cm]
 $0.4$ &  $2.07032004603564$ & $4.82\times10^{-4}$  & $2.43\times10^{-7}$ & $5.83\times10^{-8}$   \\[0,15cm]
 $0.5$ &  $2.10653065971263$ & $5.45\times10^{-4}$  & $2.75\times10^{-7}$ & $6.59\times10^{-8}$   \\
\bottomrule
\end{tabular}
\end{table}\bigskip
%

\subsection*{Example 2}
The unique solution to the IVP 
\begin{equation}\label{eq:ex02}
 \left\{\begin{array}{l}
 y'=y^2, \\[0,3cm] y(0)=1,
\end{array}\right. 
\end{equation}
is given by $y(t)=(1-t)^{-1}$.
In this scenario, the integral in procedure outlined by Eq.\,\eqref{eq:num_SIO} can also be calculated analytically.
The errors between the analytical solution and the numerical solutions using $T_2(t-h,t)$, $T_4(t-h,t)$, and  SIO(4) are presented in Table \ref{tab:example_02}.
%
\begin{table}[H]
\footnotesize
\centering
\parbox{12cm}{\caption{Absolute error between analytical and numerical solutions for IVP $y'=y^2$, $y(0)=1$, using $T_2(t-h,t)$, $T_4(t-h,t)$, and  SIO(4) ($m=3$), with $h=0.1$.}\label{tab:example_02}}
\begin{tabular}{ccccc}
\toprule
\multicolumn{2}{c}{Solution} &\multicolumn{3}{c}{Error} \\
\cmidrule(r){1-2} \cmidrule(r){3-5}
  $t$ & $y(t)$  & $T_2(t-h,t)$ & $T_4(t-h,t)$ & SIO(4) \\
 \midrule
 $0.1$ &  $1.11111111111111$ & $1.11\times10^{-3}$  & $1.11\times10^{-5}$ & $1.35\times10^{-6}$   \\[0,15cm]
 $0.2$ &  $1.25000000000000$ & $3.11\times10^{-3}$  & $3.52\times10^{-5}$ & $4.34\times10^{-6}$   \\[0,15cm]
 $0.3$ &  $1.42857142857143$ & $6.82\times10^{-3}$  & $8.96\times10^{-5}$ & $1.13\times10^{-5}$   \\[0,15cm]
 $0.4$ &  $1.66666666666667$ & $1.40\times10^{-2}$  & $2.21\times10^{-4}$ & $2.85\times10^{-5}$   \\[0,15cm]
 $0.5$ &  $2.00000000000000$ & $2.91\times10^{-2}$  & $5.75\times10^{-4}$ & $7.72\times10^{-5}$   \\
\bottomrule
\end{tabular}
\end{table}\bigskip
%

The results presented in Section \ref{sec:SIO_order} demonstrated that the order depends only to the calculation of the integral in the Eq.\,\eqref{eq:num_SIO}. 
In the following example, the integral will be computed using Gauss Quadrature \cite{burden_2015} with $m=3$, corresponding to order 5. However, other methods with orders greater than or equal to 5 yielded the same results, given the precision of decimal places used.
%

%
\subsection*{Example 3}
Let the IVP 
\begin{equation}\label{eq:ex03}
 \left\{\begin{array}{l}
 y'=\dfrac{1}{3y^2}, \\[0,3cm] y(0)=1.
\end{array}\right. 
\end{equation}\bigskip
The unique solution of PVI \eqref{eq:ex03} is $y(t)=(t+1)^{1/3}$.
In this case, the integral in the procedure given by Eq.\,\eqref{eq:num_SIO} cannot be calculated analytically in terms of elementary expression. However, a numerical scheme can be employed to calculate the integral and obtain approximations of $w_{i+1}^{[\ell]}$ in Eq.\,\eqref{eq:num_SIO}, such as employing Gauss-Quadrature of order 5. 
The errors between the analytical solution and the numerical solutions using $T_2(t-h,t)$, $T_4(t-h,t)$, and  SIO(4) are shown in Table 
\ref{tab:example_03}.
%
%
\begin{table}[H]
\footnotesize
\centering
\parbox{12cm}{\caption{Absolute error between analytical and numerical solutions for IVP $y'=1/(3y^2)$, $y(0)=1$, using $T_2(t-h,t)$, $T_4(t-h,t)$, and  SIO(4) ($m=3$), with $h=0.1$.}\label{tab:example_03}}
\begin{tabular}{ccccc}
\toprule
\multicolumn{2}{c}{Solution} &\multicolumn{3}{c}{Error} \\
\cmidrule(r){1-2} \cmidrule(r){3-5}
  $t$ & $y(t)$  & $T_2(t-h,t)$ & $T_4(t-h,t)$ & SIO(4) \\
 \midrule
 $0.1$ &  $1.03228011545637$ & $3.84\times10^{-6}$  & $2.17\times10^{-8}$ & $1.12\times10^{-8}$   \\[0,15cm]
 $0.2$ &  $1.06265856918261$ & $6.34\times10^{-6}$  & $3.32\times10^{-8}$ & $1.79\times10^{-8}$   \\[0,15cm]
 $0.3$ &  $1.09139288306111$ & $8.00\times10^{-6}$  & $3.93\times10^{-8}$ & $2.20\times10^{-8}$   \\[0,15cm]
 $0.4$ &  $1.11868894208140$ & $9.10\times10^{-6}$  & $4.24\times10^{-8}$ & $2.44\times10^{-8}$   \\[0,15cm]
 $0.5$ &  $1.14471424255333$ & $9.83\times10^{-6}$  & $4.38\times10^{-8}$ & $2.57\times10^{-8}$   \\
\bottomrule
\end{tabular}
\end{table}\medskip
%

%
Although the results in Examples 1 and 2 were obtained from the explicitly integration in the Eq.\eqref{eq:num_SIO}, the same values can be obtained using the numerical scheme employed in Example 3.
It can be observed that the errors obtained with SIO(4), which utilizes a Taylor polynomial of order 2, are smaller than the errors obtained with Taylor's method of order 4. 
%

%
%
\section{Conclusions}\label{sec:conclusions}

The central aspect of this work lies in the utilization of the spline function $S_m(t_i,t,w)$, describe in Eq.\,\eqref{eq:spline}, where $w$ is obtained  as an approximation to the analytical solution $y(t_{i+1})$ via the iterative recursive method outlined in Eq.\,\eqref{eq:num_SIO}, employing the proposed Spline-Integral Operator (SIO) \eqref{eq:G_operator}.
The SIO($m+1$) method \eqref{eq:G_opt_ti} is an implicit method that only requires information about the initial condition to be initialized. In contrast, other implicit methods of higher order, such as Predictor-Corrector methods, need a set of initial conditions for initialization, typically obtained using other numerical methods from the initial condition. 
Nevertheless, this method relies on computing of $m\!-\!1$ derivatives of the theoretical solution $y(t)$ to obtain the $S_m(t_i,t,w)$ spline.
These derivatives of the theoretical solution can be obtained from $f(t,y)$ by applying the RHS of Eq.\,\eqref{eq:taylor_m-1}, similar to how Taylor methods handle the numerical solution of IVPs.
Moreover, determining the fixed point of SIO($m+1$) necessitates computing the integral stated in  Eq.\,\eqref{eq:G_operator}.
However, straightforward methods such as Gauss Quadrature, as utilized here, or even Simpson's methods, can be employed to compute the integral of $f(t,S_m(t_i,t,w))$ over each interval $[t_i,t_i+h]$ without compromising the method's order.
As a significant aspect of this work, we demonstrated in Proposition \ref{prop:G_h(w)_approx_m} that SIO($m+1$) method is an approximation of order $m$. Furthermore, in Corollary \ref{coro:G_h(w)_approx_m}, we established that it becomes an approximation of order $m+1$ when the fixed point is reached.
Additionally, through examples, we illustrated  that SIO($m+1$) method achieves better approximations while requiring fewer derivative calculations of the theoretical solution compared to Taylor's methods of the same order.
%

\bibliographystyle{unsrt}
\bibliography{references.bib}
%
\end{document}